\documentclass[]{article}

\usepackage{amsmath}
\usepackage{amssymb}
\usepackage{amsfonts}
\usepackage{time}
\usepackage{cases}
\usepackage{authblk}

\usepackage{epstopdf}
\usepackage{subfigure}

\usepackage[numbers,sort&compress]{natbib}
\bibpunct[, ]{[}{]}{,}{n}{,}{,}

\pdfoutput=1

\newtheorem{theorem}{Theorem}[section]
\newtheorem{lemma}[theorem]{Lemma}
\newtheorem{corollary}[theorem]{Corollary}
\newtheorem{proposition}[theorem]{Proposition}

\newtheorem{definition}[theorem]{Definition}

\usepackage{graphicx}
\usepackage{graphics}
\usepackage{subfigure}
\usepackage{epsfig}

\def\QED{\mbox{\rule[0pt]{1.5ex}{1.5ex}}}
\def\proof{\noindent{\bf Proof: }}
\def\endproof{\hspace*{\fill}~\QED\par\endtrivlist\unskip}

\DeclareGraphicsExtensions{.jpg,.png,.pdf,.ps}
\graphicspath{{./img/jpg/},{./img/png/},{./img/pdf/},{./img/svg/},{./img/ps/}}

\begin{document}

\title{An alternative approach for stability analysis  of discrete time nonlinear dynamical systems}

\author{
{R. Bouyekhf\textsuperscript{a}$^{\ast}$\thanks{$^\ast$CONTACT R. Bouyekhf. Email: rachid.bouyekhf@utbm.fr}
and L.T. Gruyitch\textsuperscript{b}}\\
{\textsuperscript{a}Nanomedicine Lab, Universit\'e Bourgogne Franche-Comt\'e, UTBM, F-90010 Belfort, France\\
\textsuperscript{b}John Kennedy 31/15, Belgrade 11070, Serbia}
}
\date{}
\maketitle

\begin{abstract}

The paper endeavours to solve the problem of the  necessary and sufficient conditions for  testing asymptotic
stability of the equilibrium state  without using a positive definite or semi-definite Lyapunov
function  for time-invariant nonlinear discrete-time dynamical systems. The
solution is based on the concept of the  $\bf G$-functions  introduced in this paper. As a result,
 new  necessary and sufficient  conditions for asymptotic stability of such systems
and an estimation or the exact determination   of the asymptotic  stability domain of the
 state $x=0$ are established. 
Examples are worked out to illustrate the results.
\end{abstract}

\paragraph{keywords}
Discrete-time nonlinear systems, asymptotic stability,
stability domain, $\bf G$-functions.

\section{Introduction}

So far, the study of the stability of discrete-time nonlinear systems has been elaborated by the Lyapunov method by following and accepting its basis established for continuous-time systems, which is the concept of  definite \cite{Kalman} and semi-definite  functions \cite{Iggidr}, \cite{LaSalle1}.
It is well known that two of the main properties of such functions
are their continuity and their global minimum (usually zero) value at
the equilibrium state. Hence, the Lyapunov method reduces the stability analysis to a search for a  Lyapunov function which is decreasing  along trajectories of the system as long as $x\neq0$. Unfortunately, even in some simple cases, it can be very difficult to find such an appropriate  function.

With the best of our acknowledge,  the efforts that have been made up to the present time in the framework of discrete-time systems are limited only to some special classes of systems (linear, polynomial systems \cite{Giesl}, \cite{Oschea}, homogeneous \cite{Rosier},  exponential stable systems \cite{Peter}, \cite{Giesl1}).   
Generally speaking, for a given discrete time nonlinear system  one has little hope to find the Lyapunov function unless the system belongs to one of the special classes. 
Hence,  two questions arise naturally: When the method of the first approximation fails, why the stability test must be restricted in general to only Lyapunov  functions? Is it possible to broaden Lyapunov's  method    to certain functions without requiring their positive definiteness or semi-definiteness in the framework of discrete-time nonlinear systems? Such a possibility will enable us to enlarge the scope of the stability to other functions which have no mathematical reason to be excluded.

This paper endeavours to give  positive  answer
to  the above mentioned question. More precisely, we present analytical  tests for the asymptotic stability and an estimation or the exact determination of the asymptotic stability domain of the equilibrium state  for discrete-time nonlinear systems. The tests have to be carried out by using the concept of the $\bf G$-functions rather than that of the positive
definite or semi-definite Lyapunov functions.

The paper is organized as follows. Section 2 is devoted to the terminology and notations and summarizes some  results that are needed in the subsequent development. Section 3 defines the concept of $\bf G$-functions on which our approach is based and discusses their properties.  The final results are given in Section 4 with
 illustrative examples.

\section{System description and definitions}

We consider discrete-time nonlinear systems of the form
\begin{equation}\label{eq1}
x(k+1)=f(x(k)),\quad x(0)=x_0,
\end{equation}
where $x(k)\in {\mathbb R}^n$ is the state vector  at
discrete time $k\in {\mathbb N}$ defined on an open neighbourhood $S\subseteq \mathbb{ R}^n$ of $x=0$ and $f:S\rightarrow
S$  is a  nonlinear vector function such that $f(0)=0$. 
The
solution  of the system \eqref{eq1}, which starts at $ x_0\in S$ is defined by $x(k ;x_0)\equiv x(k)=f^k(x_0)$, where
$f^k$ is the $k$-th multiple composition of the mapping $f$.
Throughout the paper, we shall use  the following notations: $\mathbb{R}_{\geq0}=[0,+\infty)$ and $\mathbb{R}_{>0}=(0,+\infty)$. The set
$B_{\varepsilon }=\{x:\Vert x\Vert <\varepsilon\}$  is an open ball with the radius $\varepsilon$ and center at the origin. A neighbourhood of
$x=0$ is a set $N$, which contains an open ball $B_\epsilon$, i.e.; $0\in B_\epsilon\subset N$. We denote by
$C\ell A,\;\partial A$  the closure  and the boundary of the set $A$, respectively. 
Let $V:\mathbb{R}^n \longrightarrow \mathbb{R}$  be a real  function. $\Delta V(x)=V(f(x))-V(x)$ denotes the variation of the function $V$ along the solution of the system. 
For the sake of clarity we state the following:

\begin{definition}\label{def1}{\rm
 The equilibrium  state $x=0$ of system \eqref{eq1} is:
 \begin{enumerate}
\item[{(i)}]
 stable  if, and only if, for  every
                 $\varepsilon >0$ there exists
                $\delta
                (\varepsilon )>0$ such that $\Vert x_0\Vert
                <\delta
                (\varepsilon )$ implies
    $
        \Vert x(k;x_0)\Vert <\varepsilon ,\forall k\in {\mathbb
        N}
 $.

\item[{(ii)}]
attractive  if, and only if,   their exists
                 $\eta >0$ and for every $\rho >0$
                 there exists
                $K
                (x_0,\rho )>0$ such that $\Vert x_0\Vert
                <\eta
                $ implies
   $
        \Vert x(k;x_0)\Vert <\rho ,\hskip 4mm\forall\;
        k>K(x_0,\rho),
 $
or equivalently: $\displaystyle
\lim_{k\rightarrow+\infty}\|x(k,x_0)\|=0$.
\item[{(iii)}] asymptotically stable if, and only if, it is both stable and
attractive.
\end{enumerate}
}
\end{definition}

The notions of the domain of stability, the domain of attraction and the domain
of asymptotic stability used in this
paper preserve the same meaning as proposed and used in
 \cite{Grujic}, \cite{Hahn}, \cite{Miller}.

\begin{definition}\label{def2}{\rm
 The equilibrium  state $x=0$ of system \eqref{eq1} has
the domain of:
 \begin{enumerate}
\item[1)] stability $ D_s\subseteq  S$ if, and only if, both:
\begin{itemize}
  \item[(i)] for every $\varepsilon >0
                      $ there is a neighbourhood  $D_s(\varepsilon)\subseteq S$ of $x=0$ depending on $\varepsilon$
                      such that   $\Vert x(k;x_0)\Vert <\varepsilon$
                      for all $k\in {\mathbb N}$,
                       provided only that $x_0\in D_s(\varepsilon )$,
  \item[(ii)] ${ D}_s=\displaystyle\bigcup_{\varepsilon
                       >0}D_s(\varepsilon)$.
\end{itemize}

\item[2)] attraction $D_a\subseteq S$ if, and only if, the set
$D_{a}$ is the  largest neighbourhood of $x=0$ such that for every $%
x_{0}\in D_{a}$ and $\rho>0$  there exists $%
K\left( x_{0},\rho \right) >0$ guaranteeing that $\left\Vert x\left(
k;x_{0}\right) \right\Vert <\rho $ for every $k>$ $K\left( x_{0},\rho
\right)$.

\item[3)] asymptotic
stability $D_{as}\subseteq S$ if, and only if, both:
\begin{itemize}
  \item[(i)] $D_{as}=D_s\cap D_a$,
  \item[(ii)] the set $D_{as}$ is a neighbourhood of $x=0$.
\end{itemize}
\end{enumerate}
}
\end{definition}

\paragraph{Remark 1:} 
It is important to note that, the domain of stability $D_s$ introduced in Definition \ref{def2} is defined as the union of all $D_s(\varepsilon)$, where $D_s(\varepsilon)$ is  a neighbourhood  of $x = 0$ relative to each $\varepsilon>0$. Besides, the domain of attraction is defined as the  largest neighbourhood of $x=0$ that satisfies condition 2).

Notice  that, the invariance of the asymptotic stability domain is an essential element of the stability concept. Three propositions given here will be employed in the next sections. Although the results presented seem to be classical in nature  (see our earlier work on the subject \citep{Rachid}), proofs are included here in order to render the paper self-contained.

\begin{proposition}\label{prop1}
The stability domain $D_s$ of the state $x=0$ of the system \eqref{eq1} is
positively invariant with respect to system motions, that is, $ x_0\in  D_s$
implies $x(k;x_0)\in  D_s$ for every $k\in {\mathbb N}$.
\end{proposition}

\proof
 Let $\varepsilon >0$ and $x_0\in D_s(\varepsilon )$.  If the
proposition were not true, there would exist an instant $k_b\in {\mathbb
N}$ for which  $x(k_b)=x(k_b;x_0)\notin D_s(\varepsilon )$. Hence,
$\exists\; k_c\geq k_b$ such that $\Vert x(k_c;x(k_b))\Vert \geq
\varepsilon $. But $ x(k_c;x(k_b)) = x(k_c;x_0)$. Therefore,
$\|x(k_c;x_0)\| \geq \varepsilon $, which implies $x_0\notin
D_s(\varepsilon )$ that contradicts $x_0\in D_s(\varepsilon )$ and
proves the proposition.\endproof

\begin{proposition}\label{prop2}
 The domain of attraction $D_a$ of the state $x=0$ of the system \eqref{eq1}
  is positively invariant with respect to system motions.
\end{proposition}

\proof
 Let  $x_0\in D_a$,  then $\displaystyle
\lim_{k\rightarrow+\infty}\|x(k,x_0)\|=0$. Since $
x(k+j;x_0)=f^{k+j}(x_0)=f^{k}(f^j(x_0))=f^{k}(x(j))=x(k;x(j;x_0))$ for all $j<k$, it follows
$\displaystyle
\lim_{k\rightarrow+\infty}\|x(k,x(j;x_0))\|=0$. Hence,
$x(j;x_0)\in D_a$ for all $j<+\infty$, which proves the
proposition.\endproof

\begin{proposition}\label{prop3}
 The domain of asymptotic stability $D_{as}$ of the state $x=0$ of the
system \eqref{eq1} is positively invariant with respect to system motions.
\end{proposition}

Proposition \ref{prop3}  is a direct consequence of Definition \ref{def2}, Proposition  \ref{prop1}  and Proposition  \ref{prop2}.

\begin{definition}\label{def3}{\rm
 A set $E$ is an  estimate of the asymptotic stability domain $D_{as}$ if, and only if:

\begin{enumerate}
\item[{(i)}]  $E$ is a neighbourhood of $x=0$,

\item[{(ii)}]   $E\subseteq D_{as}$ and,

\item[{(iii)}]  $E$ is positively invariant set of the system.
\end{enumerate}
}
\end{definition}

\section{Definition of the $\bf G$-functions}\label{G-f}
We introduce here the class of $\bf G$-functions (not to be confused with the Meijer or Siegel $\bf G$-functions  \cite{Luke}) 
that will serve to stability analysis without requiring their positive definiteness or semi-definiteness. Let
   the extended real valued function $g : {\mathbb R}^n\rightarrow \widehat{{\mathbb R}} = {\mathbb R}\cup\{-\infty\} $ be given. Needless to say that, such a function may  be also real valued, but the  extended real valued function used here simply allows   the limit operation to take  $-\infty$ into account. Now, denote by $\mathcal{R}(g)$ the range of the function $g$ and  for   $\zeta\in \mathcal{R}(g)$, denote 
$
G_\zeta=\{x\in {\mathbb R}^n: g(x) < \zeta \}
$
 the $\zeta$-sub-level set of $g$.

\begin{definition}\label{propG:} {\rm We say  that, the  extended real valued function $g : {\mathbb R}^n\rightarrow \widehat{{\mathbb R}}$ is a $\bf G$-function if, and only if, it obeys:
\begin{enumerate}
\item[1)]  Either, $\displaystyle\min_{x\in {\mathbb R}^n}g(x)=g(0)$ if $g(0)$ is defined or 
$\displaystyle\lim_{x \rightarrow 0}g(x)=-\infty$,
\item[2)]  There is no point $x^*\in {\mathbb R}^n$ such that  $\displaystyle\lim_{x\rightarrow x^*}g(x)=+\infty$,
\item[3)] $g(x)\leq \zeta$ if, and only if, $x\in C\ell G_\zeta$  for all
$ \zeta >\zeta_m\in \mathcal{R}(g)$, where $\zeta_m=g(0)$
if $g(0)$ is defined or $\zeta_m=\displaystyle\lim_{x \rightarrow 0}g(x)=-\infty$,
\item[4)]  $g(x)=\zeta$ if, and only if, $x\in \partial G_\zeta$ for all
$ \zeta>\zeta_m\in \mathcal{R}(g)$.
\end{enumerate}}
\end{definition}

\paragraph{Remark 2:}  It is important to note that, no   continuity, positive definiteness or semi-definiteness conditions  are required for the  $\bf G$-functions. This explains that the stability conditions expressed in terms of a $\bf G$-function are substantially  different than the existing Lyapunov stability criteria for discrete-time nonlinear systems.

Now,  in order to show the meaning of the ${\bf G}$-functions,  let us present some examples. Let the following extended real valued functions $g_1,g_2 : {\mathbb R}^n\rightarrow \widehat{{\mathbb R}} $ be
 in the form $$g_1(x)=
\begin{cases} \frac{\ln(\|x\|)}{\|x\|}& \mbox{if } x\neq0\in {\mathbb R}^n, \\
-\infty&\mbox{if } x=0.\\ 
\end{cases},\quad g_2(x)=
\begin{cases} \ln(\|x\|)& \mbox{if } x\neq0\in {\mathbb R}^n, \\
-\infty&\mbox{if } x=0.\\ 
\end{cases}$$
They are discontinuous at $x=0$ and neither  positive definite nor semi-definite.  These functions are the 
${\bf G}$-functions 
 because in both cases $\zeta _m=-\infty$. Obviously, 
 conditions 2),  3) and 4) are  verified.

Another  examples where the ${\bf G}$-function is continuous everywhere are 
$$
g_1(x)=
\begin{cases} -\frac{\sin(\|x\|)}{\|x\|}& \mbox{if } x\neq0\in {\mathbb R}^n,\\
-1&\mbox{if } x=0. 
\end{cases}\quad g_2(x)=-\cos({\scriptstyle{\sum_i}} x_i)\exp(-\|x\|^2).
$$
These functions are not positive definite or semi-definite. Both  have a global  minimum at $0$, so that $\zeta_m=g(0)=-1$.  Clearly, the other conditions  are  satisfied.

Finally, let us present  an example where the ${\bf G}$-function is discontinuous at some $x$. Indeed,
let the following extended real valued function $g : {\mathbb R}^n\rightarrow \widehat{{\mathbb R}} $ 
$$g(x)=
\begin{cases} \ln(\|x\|)& \mbox{if } \|x\|\neq0<1 \\
-\infty&\mbox{if } x=0,\\ 
\|x\|& \mbox{if } \|x\|\geq 1.
\end{cases}$$
Clearly, this function is not continuous at the set $ A=\{x\in {\mathbb R}^n: \|x\|= 1\}\cup \{0\}$.  Of course, it is a ${\bf G}$-function because all conditions  of Definition  \ref{propG:} are  fulfilled.

In summary, we note that the class of $\bf G$-functions requires only  either a global minimum at zero (not necessary equal to 0)  or     
$\displaystyle\lim_{x \rightarrow 0}g(x)=-\infty$.  Conditions 3) and 4) imply simply that 
$ G_{\zeta _1}\subset G_{\zeta _2}$ for  all $\zeta _1<\zeta _2$ in $\mathcal{R}(g)$. 
On the other hand, condition 2) ensures, among other things, the connectedness of the set $ G_{\zeta}$  because the connectedness is an essential element for the proof of the stability. It is not restrictive because there is a large range of functions that satisfy this condition.

\section{Main results}

We begin this section by stating some technical lemmas necessary for the proofs of
the main results of this paper.

\begin{lemma} \label{lem1}
 Let $g$ be a  $\bf G$-function. In order that the  $\zeta$-sub-level set $G_\zeta$  to be positively invariant with respect to
motions of the system \eqref{eq1} it is  sufficient that
 the function $g$
obeys
\begin{equation}\label{eqbis3}
\Delta g(x)=g(f(x))-g(x)\leq 0,\hskip 5mm\forall x\neq0 \in G_\zeta.
 \end{equation}
\end{lemma}
\proof
Let $%
x_0\in G_\zeta$  be arbitrary. Suppose that $%
\exists k_i\in \mathbb{N} $  such that $x_i=x(k_i;x_0)\in G_\zeta$
and $x_{i+1}=x(k_{i}+1;x(k_i))=f(x_i)\not \in G_\zeta$. We then have $g (x_i)<\zeta $ and $g
(f(x_i))\geq \zeta $. Hence,
$
\Delta g (x_i)=g(f(x_i))-g (x_i)>0.
$
This result contradicts condition \eqref{eqbis3} of the Lemma and 
completes the proof.\endproof

\begin{lemma} \label{lem2}
 Let $g$ be a  $\bf G$-function and  let  $G^c_\zeta={\mathbb R}^n\backslash G_\zeta
 $ be the complement in
${\mathbb R}^n $ of $G_\zeta$. If the   function $g$
obeys
\begin{equation}\label{eqbis12}
\Delta g(x)\geq 0,\hskip 2mm\forall x \in G^c_\zeta,
\end{equation}
then $G^c_\zeta$ is positively invariant set with respect to  motions
of system \eqref{eq1}.
\end{lemma} 
\proof The proof follows the analogous  arguments  to those of the proof of Lemma \ref{lem1}.\endproof

We are now ready to state our stability criterion. To this end, recall that $\zeta_m$ introduced in Definition \ref{propG:} is defined as $\zeta_m=g(0)$
if $g(0)$ is defined or $\zeta_m=\displaystyle\lim_{x \rightarrow 0}g(x)=-\infty$.
 
\begin{theorem}\label{thm1} 
The existence  of   $\bf G$-function $g$ and  a number 
 $\Lambda \in (\zeta _m,+\infty)$ such that $G_\Lambda\subseteq S$ is an open connected neighbourhood of $x=0$ and
\begin{equation}\label{diff}
\Delta g(x)=g(f(x))-g(x)<0,\hskip 2mm\forall x\neq 0\in G_\Lambda,
\end{equation}
 are  necessary and sufficient conditions for the state $x=0$ of
the system \eqref{eq1} to be  asymptotically stable and $G_\Lambda$ to be an estimate of
its domain of the asymptotic stability $D_{as}$, i.e., $G_\Lambda\subseteq D_{as}$. 
\end{theorem}

\proof 
Necessity: Let $x=0$ of system \eqref{eq1} be asymptotically stable. From converse theorem of Lyapunov (see \citep{Lakshmikantham}),  we know that there exists a positive definite function $V$ defined on a neighbourhood $\mathbb{U}$ of $x=0$ 
such that along the trajectories of the system  we have $\Delta V(x)=V(f(x))-V(x)<0,\;\forall x\neq 0\in D_{as}$. Besides, the  $\alpha$-sub-level set of $V$, $\mathbb{V}_\alpha=\{x\in \mathbb{U}:  V(x)<\alpha\}$
 is a subset of $ D_{as}$ and  positively invariant for some $\alpha >0$.
From this information, we have to show the existence of a number $\Lambda$ such that the ${\bf G}$-function  to be proposed gives rise to an open connected set $G_\Lambda$. Indeed, note first that
  the set   $\mathbb{V}_\alpha$ might not be connected in general  and the  connectedness of $\mathbb{V}_\alpha$ depends on the nature of $V$. 
Let   $\sigma $ be such that $0<\sigma\leq\alpha$ and the $\sigma$-sub-level set 
$\mathbb{V}_\sigma=\{x\in \mathbb{U}:  V(x)<\sigma\}$ is the largest open connected subset of $\mathbb{U}$ containing $x=0$. Because $V$ is positive definite, such a number $\sigma$ exists since $V$ is continuous and $x=0$ is the global minimum of $V$.
It is clear that    $\mathbb{V}_\sigma\subseteq \mathbb{V}_\alpha \subseteq D_{as}$. Furthermore, $\mathbb{V}_\sigma$ is positively invariant with respect to motions of the system.  With this in mind, let 
now  the extended real valued function $g : {\mathbb R}^n\rightarrow \widehat{{\mathbb R}} $  be
 in the form
 $$g(x) =\begin{cases} \ln (V(x))  &\mbox{if } x\neq0\in \mathbb{V}_\sigma,\\-\infty& \mbox{if } x=0,\\
\ln(\sigma)& \mbox{if } x\notin \mathbb{V}_\sigma.
\end{cases}$$
It  is neither continuous at $x=0$ 
nor positive definite or semi-definite. We shall first prove that this function is a ${\bf G}$-function. Indeed, 
 observe that $V(0)=0$ implies $\displaystyle\zeta_m=\lim_{x \rightarrow 0}g(x)=-\infty$,
so that the condition 1) of Definition \ref{propG:} is verified. On the other hand, since  for $x\in \mathbb{V}_\sigma$ we have
$\ln (V(x))<\ln(\sigma)<\sigma<+\infty$, then the
 condition 2) of Definition \ref{propG:}  is also verified. Finally,  it is clear from the definition of $g$ that conditions 3) and 4) of Definition \ref{propG:} are trivially fulfilled.  This proves the
existence of the  ${\bf G}$-function $g$.

Next, observe that  if we take $\Lambda=\ln(\sigma)$, then 
$G_\Lambda=\{x\in \mathbb{R}^n: \quad\ln(V(x))<\Lambda\}=\mathbb{V}_\sigma $.
Hence,  $G_\Lambda$ is an open  connected neighbourhood of $x=0$.  
This proves the
existence of $G_\Lambda \subseteq D_{as}$ and  $G_\Lambda \subseteq S$ due to $D_{as} \subseteq S$ in view of Definition \ref{def2}.
Finally, since $\Delta V(x)<0$ for all $ x\neq 0\in \mathbb{V}_\sigma$ along the trajectories of the system  and    $\mathbb{V}_\sigma= G_\Lambda$ is positively invariant, it follows
$$V(f(x))< V(x),\quad \forall x\neq 0\in G_\Lambda.$$
 Thus, for all $x\neq 0\in G_\Lambda$ we have
$$\Delta g(x)=\ln\bigl(V(f(x))\bigl) - \ln\bigl( V(x)\bigl)  <0,$$
which proves the necessity.

Sufficiency: The hypotheses of the theorem   prove that for all $\zeta \in
(\zeta _m,\Lambda ]$ the set $G_\zeta $ is positively invariant with
respect to the motions of the system \eqref{eq1} (Lemma \ref{lem1}), i.e.,
\begin{equation}\label{eq3}
x_0\in G_\zeta \Rightarrow x(k;x_0)\in G_\zeta,\quad
\forall k\in {\mathbb N},\;\forall \zeta \in (\zeta _m,\Lambda ].
\end{equation}
Let  $\varepsilon >0$ be arbitrarily chosen and let  $\zeta
(\varepsilon )\in (\zeta _m,\Lambda ]$ be such that
$ C\ell G_{\zeta (\varepsilon )}\subset
B_\varepsilon $.
The existence of the number $\zeta(\varepsilon ) $ obeying the last
condition is assured by
 the definitions of $g$
and   $G_{\zeta(\varepsilon )} $. Now, let $\delta (\varepsilon )>0$
be such that
\begin{equation}\label{eq4}
B_{\delta (\varepsilon )}\subset C\ell G_{\zeta (\varepsilon
)}\subset B_\varepsilon. 
\end{equation}
Such $\delta (\varepsilon )$ exists because $G_{\Lambda}$ is an open connected neighbourhood of $x=0$. From
\eqref{eq3} and \eqref{eq4} we have
\begin{equation}\label{eq6}
x_0\in B_{\delta (\varepsilon )}\Rightarrow x(k;x_0)\in G_{\zeta
(\varepsilon )}\subset B_\varepsilon ,\hskip 2mm\forall k\in {\mathbb
N}.
\end{equation}
which proves the stability of $x=0$  and
$G_{\zeta (\varepsilon )}\subseteq D_s$ for all $\zeta (\varepsilon) \in (\zeta_m,\Lambda]$. Moreover, \eqref{eq6} and $\zeta
(\varepsilon) \in (\zeta_m,\Lambda]$ show that $G_\Lambda$ is
positively invariant with respect to the system motions. Hence,
$G_\Lambda$ is also an estimate of $D_s$ in view of Definition \ref{def3},
i.e., $G_\Lambda\subseteq D_s$.

In order to prove attraction we suppose opposite, i.e., that there exist
$ x_0^*\neq0\in G_\Lambda$ and $\alpha >0$ such that
$
\displaystyle\lim_{k\rightarrow +\infty }\Vert x(k;x_0^*)\Vert =\alpha
$.  
Let $\displaystyle \lim_{k\rightarrow +\infty }g(x(k;x_0^*))=\gamma$. Since $x(k;x_0^*)\neq 0$ and $\Delta g(x)< 0$ for all $x\neq0\in G_\Lambda$, it follows from 1) of
Definition \ref{propG:}  that $\gamma>  \zeta _m$ (otherwise $\gamma\leq\zeta _m$ would imply $x(+\infty;x_0^*)=0$ in view of 1)  of
 Definition \ref{propG:}, which is a contradiction since $\Vert x(+\infty;x_0^*) \Vert \neq0$ by assumption).  Now, observe that since $g$ is decreasing and $\displaystyle \lim_{k\rightarrow +\infty }g(x(k;x_0^*))=\gamma$, then  the trajectory $x(k;x_0^*)$ lies outside $G_\gamma$. It follows  for all $k\in \mathbb N$
 $$x(k,x_0^*)\in G(x_0^*,\gamma)=\{x\in \mathbb{R}^n: \gamma \leq g(x)\leq g(x_0^*)\}\subset G_\Lambda.$$
Let $-\mu=\displaystyle \sup_{x\in  G(x_0^*,\gamma)} \Delta g(x)  $. The supremum clearly exists due to  stability of $x=0$ of system \eqref{eq1} and  positive invariance of  $G_\Lambda$. Thus, $-\mu <0$ by \eqref{diff}.  Hence,
$\Delta
g(x)\leq-\mu<0$ for all $x\in G(x_0^*,\gamma)$.
This infers that for all $x(k)\in G(x_0^*,\gamma)$
$$
g(x(k;x_0^*))\leq g(x_0^*)-\sum_{i=0}^{k-1}\mu= g(x_0^*)-k\mu.
$$
Therefore, for $k\rightarrow +\infty
$ we have $g(x(k;x_0^*))\rightarrow -\infty $ due to
$\mu>0$ and $g(x_0^*)\neq +\infty$  in view of condition 2) of Definition \ref{propG:}. But $\displaystyle \lim_{k\rightarrow +\infty }g(x(k;x_0^*))=\gamma\neq-\infty$, a contradiction with the definition of $\gamma$. Consequently
$\displaystyle
\lim_{k\rightarrow +\infty }\Vert x(k;x_0)\Vert =0$ for all
$x_0\in G_\Lambda,
$
so that $G_\Lambda\subseteq D_a$ and hence $G_\Lambda\subseteq D_s\cap
D_a=D_{as}$,  which completes the proof of the theorem.
\endproof 

\paragraph{Remark 3:} As in the Lyapunov approach, the concept of $\bf G$-functions used in theorem \ref{thm1} has an intuitive geometric meaning. If a solution $x(k, x_0)$ starts at $x_0\in G_\Lambda$, then the positive invariance of $G_\Lambda$  with respect to the system ensures that the solution remains in $ G_\Lambda$ for all $k\in \mathbb{N}$. Moreover, the positive invariance of $G_\zeta$ for all $\zeta\leq\Lambda$ and  the strict inequality $\Delta g(x)<0$ will ensure also that the system solution   will progress by entering a decreasing sequence of $\zeta$-sub-level-sets $ G_{\zeta _1}\supset G_{\zeta _2}\cdots\supset G_{\zeta _i}\supset \cdots$ with $\zeta _1=g(x(0)>\zeta _2=g(x(1))>\cdots>\zeta _i=g(x(i-1))>\cdots>\zeta _m$. Consequently, in the long run the system solution will asymptotically converge to $C\ell G_{\zeta _m}=\{x\in \mathbb{R}^n: g(x)\leq \zeta _m\}=\{0\}$.

\paragraph{Remark 4:}  It is  important to note that,
if $g(0)$ is finite, one might be tempted to consider $V(x) = g(x)-g(0)$ and thus obtain a Lyapunov function. 
 However, observe in this case that $\Delta V(x)  =  g(f(x))-g(x) $. Hence, by doing so, we bring back to the fundamental problem of the analysis of the sign of $\Delta V(x) $ with respect to $g$ (which is  sign indefinite) rather than with respect to the positive definite function $V$. The routine proof of attraction of $x=0$ using $\Delta V(x) <0$ in the existing sufficient conditions  given by the Lyapunov method is based on the fact that $V$ is positive for $x\neq 0$.
 Therefore, taking $V(x) = g(x)-g(0)$  does not recapture the technique of Lyapunov function.

Next, to illustrate  the notion of the $\bf G$-functions, let us look at the following simple but illustrative examples. First of all, it is important to note  that, the asymptotic stability of  the presented examples can be tested by several tools, but the aim here is to show that the $\bf G$-functions are also applicable.

\paragraph{Example 1:}  Consider the following system  defined on $S=(-1;1)$
$$
x(k+1)=\sin (x(k))\bigl(\tanh(x(k))-x(k)\bigr)=f(x(k)).
$$
 The system has $x=0$ as a unique equilibrium state on $S$.  As shown in Section \ref{G-f}, let the following extended function 
$$g(x)=
\begin{cases} \frac{\ln(\vert x\vert)}{\vert x\vert}& \mbox{if } x\neq0\in \mathbb{R},\\
-\infty&\mbox{if } x=0.\\ 
\end{cases}$$
be candidate $\bf G$-function for such system.
Since, the function $f$ maps $S$ into itself, we  take $\Lambda=0$ so that $G_\Lambda= (-1;1)= S$ is an open connected neighbourhood of $x=0$. Now, we test   $\Delta g(x)$ on $G_\Lambda/\{0\} $. Indeed,  for all $x\neq0\in G_\Lambda$ we have
$$
\Delta g(x)=\frac{\ln(\vert f(x)\vert)}{\vert f(x)\vert}-\frac{\ln(\vert x\vert)}{\vert x\vert}
=\frac{1}{\vert f(x)\vert\vert x\vert}\Bigl(\vert x\vert\ln(\vert f(x)\vert)-\vert f(x)\vert\ln(\vert x\vert)\Bigr).$$
Let $p(x)= \vert x\vert\ln(\vert f(x)\vert)-\vert f(x)\vert\ln(\vert x\vert)$. Since $p(-x)=p(x)$, it is enough to study $p(x)$ on $(0,1)$. Indeed,  for all $x\in (0,1)$, we have
$$p(x)=x\ln(\sin(x))+x\ln(\vert\tanh(x)-x\vert)-\sin(x)\vert\tanh(x)-x\vert\ln(x).$$
But, $0<\sin(x)<1$ and $\ln(x)<0$ on $ (0,1)$, it follows
$$p(x)<x\ln(\vert\tanh(x)-x\vert)-\vert\tanh(x)-x\vert\ln(x).$$
Since   $x\neq\vert\tanh(x)-x\vert)$ on $(0,1)$, then by applying the fact that $(a-b)(\ln(b)-\ln(a))<0$ for all $a\neq b\in {\Bbb R} _{>0}$, we get after some elementary
 transformations
\begin{align*}
p(x)&<\bigl[x-\vert\tanh(x)-x\vert\bigl]\bigl[\ln(\vert\tanh(x)-x\vert)-\ln(x)\bigl]+
\ln(x)\bigl[x-2\vert\tanh(x)-x\vert\bigl]\\
&+\vert\tanh(x)-x\vert\ln(\vert\tanh(x)-x\vert).
\end{align*}
Now, remark that for  $x\in (0,1)$ we have  $\ln(\vert\tanh(x)-x\vert)<0$ and $x-2\vert\tanh(x)-x\vert>0$. We conclude that $p(x)<0$ on $ (-1,1)$ because $p(x)$ is even. Therefore,   $\Delta g(x)<0$ for all $x\neq0\in G_\Lambda$. All conditions of
Theorem \ref{thm1}  are satisfied, the state $x=0$ of the system is
asymptotically stable and the set $ G_\Lambda$ is an estimate of its
domain of asymptotic stability $D_{as}$.

\paragraph{Example 2:}   Consider the following system  defined on $S=(-\pi;\pi)$
\begin{equation}\label{eq111}
x(k+1)=\frac{\sin (x(k))}{2 \cos^3(\frac{x(k)}{3})}=f(x(k)).
\end{equation}
 The system has $x=0$ as a unique equilibrium state on $S$. As shown in Section \ref{G-f}, let the following extended function 
$$g(x)=
\begin{cases} -\frac{\sin(x)}{x}& \mbox{if } x\neq0,\\
-1&\mbox{if } x=0.\\ 
\end{cases}$$
be candidate $\bf G$-function for such system with $\zeta_m=-1$.  It is well known that the function $g$   enjoys the following properties (see \cite{Feng}, \cite{Sandor}).
\begin{itemize}
\item Jordan's inequality
\begin{equation}\label{eq11111}
-1\leq-\frac{\sin(x)}{x}\leq -\frac{2}{\pi}, \quad \text{for} \; \vert x\vert\leq\frac{\pi}{2},
\end{equation}
\item S\'{a}ndor's inequality
\begin{equation}\label{eq1111}
-\cos(\frac{x}{2})<-\frac{\sin(x)}{x}<- \frac{1+\cos(x)}{2},\quad \text{for}\;  \vert x\vert<\frac{\pi}{2}.
\end{equation}
\end{itemize}

From the 
 Jordan's inequality \eqref{eq11111},
if we  take $\Lambda=-\frac{2}{\pi}$, then $G_\Lambda= (-\frac{\pi}{2};\frac{\pi}{2})\subset S$ is an open connected neighbourhood of $x=0$. In order to test    $\Delta g(x)$ on $G_\Lambda/\{0\} $, observe that  for all $x\neq0\in G_\Lambda$ we have
$\Delta g(x)=\frac{\sin(x)}{x}- \frac{\sin(f(x))}{f(x)}.$
To exploit  inequality \eqref{eq1111}  we have to show that $f(x)\in (-\frac{\pi}{2};\frac{\pi}{2})$. To this end,   since $f(x)$ is odd, it suffices to study  $f$ on $(0, \frac{\pi}{2})$. Indeed, elementary computation yields  $f'(x)=\frac{1}{2}\cos(\frac{2x}{3} )\cos^{-4}(\frac{x}{3} )>0$ on $(0, \frac{\pi}{2})$. Therefore,  $0=f(0)<f(x)<f(\frac{\pi}{2})=\frac{4}{3\sqrt{3}}<\frac{\pi}{2}$ so that $f(x)\in (-\frac{\pi}{2};\frac{\pi}{2})$. Consequently, we can use inequalities \eqref{eq1111} to obtain
$$\Delta g(x)<\cos(\frac{x}{2})- \frac{1+\cos(f(x))}{2}.$$
But, $ \frac{1+\cos(f(x))}{2}=\cos^2\bigl(\frac{f(x)}{2}\bigl)$ and $\cos^2(\frac{x}{4})=\frac{1+\cos(\frac{x}{2})}{2}\geq \cos(\frac{x}{2})$. This implies
$$
\Delta g(x)<\cos^2\bigl(\frac{x}{4}\bigl)-\cos^2\bigl(\frac{f(x)}{2}\bigl)
=\Bigl[\cos\bigl(\frac{x}{4}\bigl)-\cos\bigl(\frac{f(x)}{2}\bigl)\Bigl]\Bigl[
\cos\bigl(\frac{x}{4}\bigl)+\cos\bigl(\frac{f(x)}{2}\bigl)\Bigl].
$$
Let now
$h(x)=\cos\bigl(\frac{x}{4}\bigl)-\cos\bigl(\frac{f(x)}{2}\bigl)$. Since, $\cos(\frac{x}{4})+\cos(\frac{f(x)}{2})>0$ on $(-\frac{\pi}{2};\frac{\pi}{2})$ and $h(-x)=h(x)$,  it is sufficient  to study $h(x)$ on $(0, \frac{\pi}{2})$. For, define 
$p(x)=\frac{x}{4}-\frac{f(x)}{2}$. After some elementary transformations
 we get
$p'(x)=\frac{1}{4}-\frac{1}{2}f'(x)=\frac{1}{4}\tan^4(x)>0.$
Therefore, $p(x) >p(0)=0$ so that $\frac{f(x)}{2}<\frac{x}{4}$ on $(0, \frac{\pi}{2})$. Since  $\cos(x)$ is decreasing on $(0, \frac{\pi}{2})$ we get $\cos(\frac{x}{4})<\cos(\frac{f(x)}{2})$  and thus $h(x)<0$ for all 
$x\in (-\frac{\pi}{2};\frac{\pi}{2})$. Consequently,  $\Delta g(x)<0$ for all $x\neq0\in G_\Lambda$.
The state $x=0$ of the system is
asymptotically stable and the set $ G_\Lambda$ is an estimate of its
domain of asymptotic stability $D_{as}$.

Next, recall that a function $\phi (\cdot ):{\Bbb R} _{\geq0}\rightarrow {\Bbb R}_{\geq0}$ is 
 a comparison function of the class ${\mathcal K}$ defined by Hahn \cite{Hahn} if, and only if, it is continuous, strictly increasing and
 $\varphi(0)=0$.

Classically, in the existing Lyapunov stability theory, it is customary to require also the existence of a function $\phi$ of class ${\mathcal K}$ such that $\Delta V(x)\leq -\phi(\| x\|)$ instead of requiring only $\Delta V(x)<0$, where $V$ is a Lyapunov function. It is worthwhile to show how this requirement can  be also extended  to $\bf G$-functions.
 In this framework, we have the next result.

\begin{theorem}\label{thm2}
The existence  of  a
$\bf G$-function $g$, a function $\phi $ of class ${\mathcal K}$  and a number $\Lambda \in (\zeta _m, +\infty)$ such that $G_\Lambda\subseteq S$ is an open connected neighbourhood of $x=0$ and 
\begin{equation}\label{difff}
\Delta
g(x)\leq -\phi(\| x\|),\hskip 2mm\forall x\neq 0\in G_\Lambda,
\end{equation}
are necessary and  sufficient conditions for the state $x=0$ of
the system \eqref{eq1} to be  asymptotically stable and $G_\Lambda$ to be an estimate of
its domain of asymptotic stability $D_{as}$, i.e., $G_\Lambda\subseteq D_{as}$. 
\end{theorem}
 
\proof By Theorem \ref{thm1}, it is clear that we need only to show  here the existence of a function $\phi $ of class ${\mathcal K}$ such that the inequality \eqref{difff} is verified. Indeed, from converse theorem of Lyapunov we know also that there exist a positive definite function $V$ and a function $\psi$ of class ${\mathcal K}$ 
such that along the trajectories of the system \eqref{eq1}  we have for all $ x\neq 0\in D_{as}$ (see \cite{Lakshmikantham},  Theorem 4.10.2)
$$\Delta V(x)=V(f(x))- V(x)\leq-\psi(\parallel x\parallel).$$
This in turn implies that  $V(x)>\psi(\parallel x\parallel), \forall x\neq 0\in D_{as}$ in view of positive definiteness of $V$. With this in mind, let  $g$ and $G_\Lambda=\mathbb{V}_\sigma$  as defined in the proof of the necessity part of Theorem \ref{thm1}. Then we have
$$V(f(x))\leq V(x)-\psi(\|x\parallel),\quad \forall x\neq 0\in G_\Lambda,$$
which implies
$$\frac{V(f(x))}{V(x)}\leq 1-\frac{\psi(\|x\|)}{V(x)},\quad \forall x\neq 0\in G_\Lambda.$$
Now, observe that since $V(x)>\psi(\|x\|)$ for all $x\neq 0\in G_\Lambda$, then $\ln\Bigl(1-\frac{\psi(\|x\|)}{V(x)}\Bigr)$ is well defined. Thus, we get
$$\Delta g(x)=\ln\Bigl(\frac{V(f(x))}{V(x)}\Bigr)\leq\ln\Bigl(1-\frac{\psi(\|x\|)}{V(x)}\Bigr)\leq -\frac{\psi(\|x\|)}{V(x)},$$
where we have used that fact that $\ln(1-x)\leq -x$ for $0\leq x< 1$. Let   $\displaystyle\gamma=\sup_{x\in G_\Lambda}V(x)>0$.  It follows
$$\Delta g(x)\leq -\frac{1}{\gamma}\psi(\|x\|),\quad \forall x\neq 0\in G_\Lambda.$$
Clearly  $\phi(\|x\|) =\frac{1}{\gamma}\psi(\|x\|)$ belongs to class $\mathcal{K}$.  This completes the proof.\endproof

Now, for the global asymptotic stability, we have the following.

\begin{corollary}
Under the conditions of Theorem \ref{thm1} or Theorem \ref{thm2}, if  $S=\mathbb{R} ^n$ and $\Lambda $ can be arbitrarily
large, then the state $x=0$ is globally asymptotically stable, i.e.,
$D_{as}=\mathbb{R} ^n$.
\end{corollary}
\paragraph{Example 3:}   Consider the following   two dimensional system defined on $\mathbb{R}^2$
\begin{align}
x(k+1)&=  x(k)y(k)\exp(-y^2(k))\\
y(k+1)&=  \alpha x(k)
\end{align}
where  $ \alpha \neq0\in \mathbb{R}$. We are interested in discovering  a condition on $ \alpha $ under which  the global asymptotic stability of the  fixed point $(0, 0)$ is guaranteed. Let the extended real valued function $g : {\mathbb R}^n\rightarrow \widehat{{\mathbb R}} $ be
 in the form 
 $$g(x)=
\begin{cases} \ln(\vert x\vert+\vert y\vert)& \mbox{if } (x,y)\neq(0,0)\in {\mathbb R}^n, \\
-\infty&\mbox{if } (x,y)=(0,0).\\ 
\end{cases}$$
 This function is the $\bf G$-function 
which evidently obeys Definition \ref{propG:}  with $\zeta
_m=-\infty$. Now, Let $\Lambda\in \mathbb{R}$ be arbitrarily chosen. Clearly,   
$
G_\Lambda=\{x\in \mathbb{R}^n :\ln(\vert x\vert+\vert y \vert)< \Lambda\},
$
is on an open connected  neighbourhood of $(x,y)=(0,0)$. We test now $\Delta g(x)$ on $G_\Lambda/\{(0,0)\}$. Indeed,
$$
\Delta g(x)
=\ln\Bigr(\frac{\vert x\vert}{\vert x\vert+\vert y \vert}\Bigr)+
\ln\Bigr(\vert y\vert  \exp(-y^2)+  \vert \alpha \vert  \Bigr).
$$
Since, $\frac{\vert x\vert}{\vert x\vert+\vert y \vert}\leq1$ and $\vert y\vert  \exp(-y^2)\leq \frac{1}{\sqrt{2e}},\; \forall y\in \mathbb{R}$, it follows
$\Delta g(x)\leq \ln\Bigr(\frac{1}{\sqrt{2e}}+\vert \alpha \vert \Bigr)$.
If $\vert \alpha \vert< 1-\frac{1}{\sqrt{2e}}$, then $\Delta g(x)<0$ for all $(x,y)\in G_\Lambda/\{(0,0)\}$. We conclude that all conditions of
Theorem \ref{thm1}  are satisfied for arbitrary $\Lambda\in \mathbb{R}$. The state $(x,y)=(0,0)$ is then globally asymptotically
stable.

So far, we have investigated the asymptotic stability problem of discrete-time nonlinear systems by means of sign indefinite $\bf G$-functions. Now, a fundamental question arises: How can the domain of asymptotic stability be completely determined by the $\bf G$-functions?  By enhancing the sufficient conditions for  asymptotic stability of the state $x=0$, the following theorem establishes a complete answer to this question.

\begin{theorem}\label{theo3}
 If all conditions of Theorem \ref{thm1} or Theorem \ref{thm2}  hold and if in addition
\begin{equation}\label{eq12}
\Delta g(x)\geq 0,\hskip 2mm\forall x\in \mathbb{R} ^n/ G_\Lambda,
\end{equation}
then the set $G_\Lambda $ is the domain of the asymptotic stability of the state
$x=0$ of the system \eqref{eq1}, i.e., $G_\Lambda\equiv D_{as}$.
\end{theorem}

\proof
By Theorem \ref{thm1} we have
\begin{equation}\label{eq13}
G_\Lambda\subseteq D_a\cap D_s=D_{as}.
\end{equation}
Let $x_0\in {\mathbb R} ^n/G_\Lambda $. Lemma \ref{lem2}  shows that ${\mathbb R} ^n/G_\Lambda $ is positively invariant, i.e., $
x_0\in {\mathbb R} ^n/G_\Lambda$ implies $ x(k;x_0)\in {\mathbb R}
^n/G_\Lambda$ for all $k\in {\mathbb N}.
$
Since ${\mathbb R} ^n/G_\Lambda $ does not contain the origin, it
follows  $\displaystyle \lim_{k\rightarrow +\infty}
\|x(k;x_0)\|\neq0$, which implies $D_a\subseteq G_\Lambda$. This and \eqref{eq13} yield $D_a\cap D_s=G_\Lambda$ so that $D_{as}= G_\Lambda$, which completes the proof.\endproof

\begin{figure}[h]
\begin{center}
  \subfigure[The graph for  $f(x)=x^2$ ]{\epsfig{figure=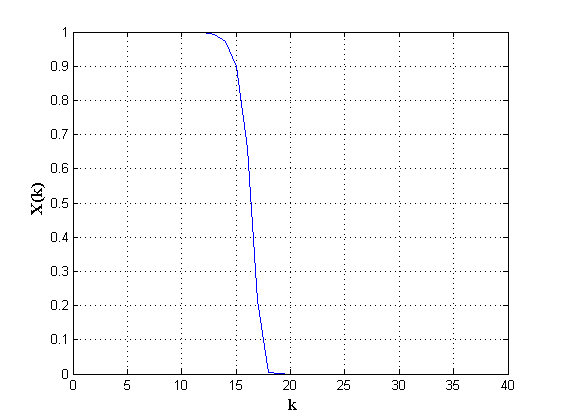,width=4.58cm}}\quad
  \subfigure[The graph for  $f(x)=e^{x^2}-e+1$ ]{\epsfig{figure=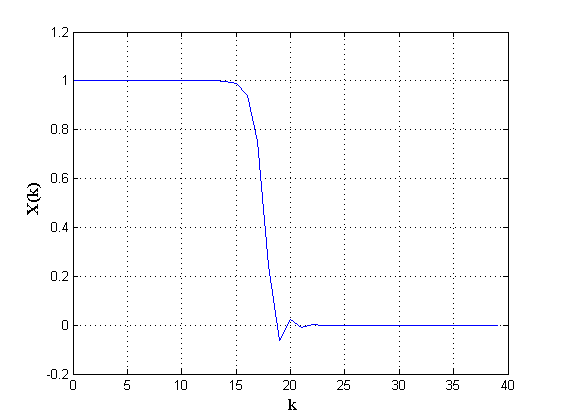,width=4.58cm}}
\end{center}
\caption{Stable trajectories}
\label{fig1}
\end{figure}

\begin{figure}[h]
\begin{center}
  \subfigure[The graph for  $f(x)=x^2$ ]{\epsfig{figure=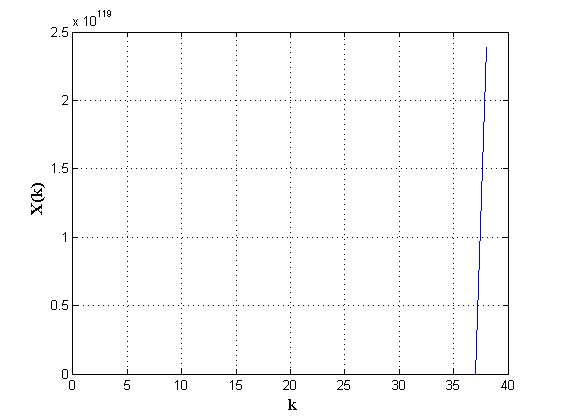,width=4.58cm}}\quad
  \subfigure[The graph for  $f(x)=e^{x^2}-e+1$ ]{\epsfig{figure=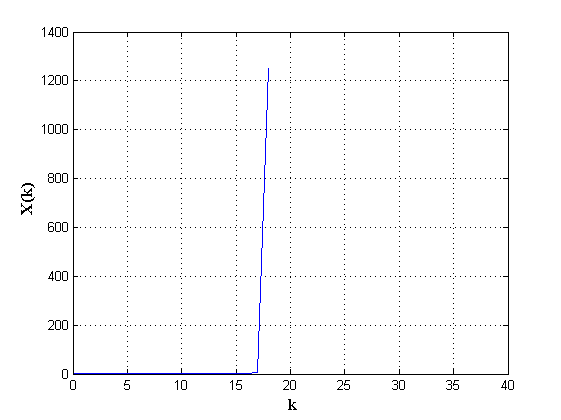,width=4.58cm}}
\end{center}
\caption{Unstable trajectories}
\label{fig2}
\end{figure}

\paragraph{Example 4:} 
Let the following family of  systems 
\begin{equation}\label{eq11}
x(k+1)=\frac{f(x(k))}{1+\vert \|x(k)\|-e^{\alpha}\vert }\;x(k),
\end{equation}
where $x(k)\in \mathbb{R}^n$,  $\alpha\in \mathbb{R}$ and $f:\mathbb{R}^n\rightarrow
\mathbb{R}$.  We are interested in discovering   functions family  $ f(x) $ under which  the set $S(\alpha)=\{x\in \mathbb{R}^n: \|x\|<e^{\alpha}\}$ is the domain of the asymptotic stability of the  fixed point $x=0$. As shown in Section \ref{G-f}, let   the following  $\bf G$-function  
$$g(x)=
\begin{cases}\ln (\| x\| )& \mbox{if } x\neq0\in {\mathbb R}^n, \\
-\infty&\mbox{if } x=0.\\ 
\end{cases}$$
be candidate  for such system.   Let $\Lambda=\alpha$, then 
$
G_\Lambda=\{x\in \mathbb{R}^n :\ln(\| x\|) < \Lambda\}=S(\alpha)
$
is an open connected neighbourhood of $x=0$.  Now, on  $G_\Lambda/\{0\}$ we have
$$
 \Delta g(x)
=\ln(\vert f(x)\vert) +\ln\Bigl(    \frac{1}{ 1+\vert \|x\|-e^{\alpha}\vert } \Bigr)
$$
According to Theorem \ref{theo3},  the set $S(\alpha)$ is  the domain of asymptotic stability if:
\begin{itemize}
\item on $S(\alpha)/\{0\} $, $\vert f(x)\vert\leq 1,$  
\item on $x\in  \partial S(\alpha)$,  $f(x)=1, $
\item on $ {\mathbb R} ^n/ S(\alpha)$, $f(x)\geq  1+ \|x\|-e^{\alpha}.$  
\end{itemize}
For the simulation purpose, consider  the case of one dimensional systems.   For $\alpha=0$, it is clear that $f(x)=x^2$ or $f(x)=e^{x^2}+1-e$ satisfy all the above conditions. In Both cases, system motions are shown in Figure \ref{fig1} for $x(0)=0.999999999\in S(0)$ and in Figure \ref{fig2} for
$x(0)=1.00000001$ out of $S(0)$. Both initial states are close to the boundary of $S(0)$. In the former
case the system motion converges to $x=0$, while in the latter case it does not. Notice also a very small difference between initial values in the two cases. The set $S(0)$
is therefore the domain of the asymptotic stability of $x=0$.



\section{Conclusion}

The theoretical study of the stability problems for the  discrete-time nonlinear
systems has been developed by using properties of a class of functions called  the $\bf G$-functions in this paper. It has resulted in the   new necessary and   sufficient
conditions for asymptotic stability of such systems and an estimation or the exact
determination of the asymptotic stability domain of the zero state
$x=0$. Hence, the paper establishes essentially an alternative approach for  
the asymptotic stability analysis  of time-invariant nonlinear
discrete-time systems. The essence is the possibility to use henceforth the sign indefinite  and discontinuous functions   for testing asymptotic stability of such systems. 

Finally,  we would like to point out that the $\bf G$-functions  approach presented here  simply offers an additional tool for the study of the stability of discrete-time nonlinear systems. It can be further developed and possibly combined with the Lyapunov function approach to generate more interesting and advanced results in the future, specifically for stabilization and control of discrete-time nonlinear systems.


\begin{thebibliography}{99}
\providecommand{\MR}{\relax\unskip\space MR }
\providecommand{\url}[1]{\normalfont{#1}}
\providecommand{\urlprefix}{Available at }

\bibitem{Rachid}
R. Bouyekhf and L.T. Grujic, \emph{Lyapunov-like solutions to stability
  problems for discrete-time systems on asymptotically contractive sets},
  Nonlinear Dynamics 18 (1999), pp. 107--127.

\bibitem{Giesl}
P. Giesl, \emph{Construction of a local and global {L}yapunov function using
  radial basis functions}, IMA J. Appl. Math 75 (2008), pp. 782--802.

\bibitem{Peter}
P. Giesl and S. Hafsteinb, \emph{Computation of {L}yapunov functions for
  nonlinear discrete time systems by linear programming}, Journal of Difference
  Equations and Applications 20 (2014), pp. 610--640.

\bibitem{Giesl1}
P. Giesl and S. Hafsteinb, \emph{Review on computational methods for {L}ypunov
  functions}, Discrete and Continuous Dynamical Systems Series B 20 (2015), pp.
  2291--2331.

\bibitem{Grujic}
L.T. Gruyitch, J.P. Richard, P. Borne, and J.C. Gentina, \emph{Stability
  Domains}, Chapman and Hall/ CRC press, London, 2003.

\bibitem{Hahn}
W. Hahn, \emph{Stability of Motion}, Springer-Verlag, Berlin, 1967.

\bibitem{Iggidr}
A. Iggidr and M. Bensoubaya, \emph{New results on the stability of
  discrete-time systems and applications to control problems}, Journal of
  Mathematical Analysis and Applications 219 (1998), pp. 392--414.

\bibitem{Kalman}
R.E. Kalman and J.E. Bertram, \emph{Control system analysis and design via the
  second method of {L}yapunov, part {II}, discrete systems}, Transactions of
  the ASME Journal of Basic Engineering 82 (1960), pp. 394--400.

\bibitem{Lakshmikantham}
V. Lakshmikantham and D. Trigiante, \emph{Theory of Difference Equations with
  Applications to Numerical Analysis}, Academic Press, San Diego, 1988.

\bibitem{LaSalle1}
J.P. LaSalle, \emph{The Stability of Dynamical Systems}, SIAM, Philadelphia,
  Pennsylvania, 1976.

\bibitem{Luke}
Y.L. Luke, \emph{The Special Functions and their Approximations}, Academic
  Press, New York, 1969.

\bibitem{Miller}
R. Miller and A. Michel, \emph{Ordinary Differetial Equation}, Academic Press,
  New York, 1982.

\bibitem{Oschea}
R. O'schea, \emph{The extension of the {Z}obov's method to sampled data control
  systems described by nonlinear autonomous difference equations}, IEEE
  Transactions on Automatic Control January (1964), pp. 62--70.

\bibitem{Feng}
B.N.G. Qi. Feng~Qi D-W.~Niu, \emph{Refinements, generalizations, and
  applications of {J}ordan's inequality and related problems}, Journal of
  Inequalities and Applications 2009, Article ID 271923 (2009), p. 52 pages.

\bibitem{Rosier}
L. Rosier, \emph{Homogeneous {L}yapunov function for homogeneous continuous
  vector fields}, Systems and Control Letters 19 (1992), pp. 467--473.

\bibitem{Sandor}
J. S\'{a}ndor, \emph{Certain trigonometric inequalities}, Octogon Mathematical
  Magazine 9 (2001), pp. 331--336.

\end{thebibliography}
\end{document}